\documentclass[11pt]{amsart}
\usepackage{epsfig,a4}

\newtheorem{theorem}{Theorem} 

\newtheorem{cor}[theorem]{\sc Corollary}
\newtheorem{prop}[theorem]{\sc Proposition}
\theoremstyle{remark}
\newtheorem{problem}{Problem}

\newenvironment{epigraph}{\begin{flushright}% 
\sl\footnotesize}%
{\end{flushright}\medskip}
\newcommand{\source}[2]{\\ \medskip --- {\sc #1\/}, #2.} 

\renewcommand{\frac}{\tfrac}
\newcommand{\vpic}[1]{{\rule[-.7em]{0em}{2em}\:\vcenter{\epsffile{#1.eps}}\:}}

\newcommand{\Z}{\mathbb{Z}}

\newcommand{\Wh}{{W\!h}}
\newcommand{\half}{\frac{1}{2}}
\renewcommand{\topfraction}{.75}

\title{On the first two Vassiliev invariants}
\author{Simon Willerton}
\email{simon@ma.hw.ac.uk}
\address{Department of Mathematics, Heriot-Watt University, Edinburgh EH14 4AS, Scotland.}

\begin{document}

\begin{abstract}
The values that the first two Vassiliev invariants take on
prime knots with up to fourteen crossings are considered.  This leads
to interesting fish-like graphs.  Several results about the values
taken on torus knots are proved.
\end{abstract}
\maketitle

\begin{epigraph}
    `First the fish must be caught.'\\
    That is easy: a baby, I think, could have caught it.
  \source{The Red Queen}{Through the Looking Glass}%
    \end{epigraph}

\section*{Introduction.}
The two simplest non-trivial Vassiliev knot invariants (see
\cite{Vas,BirmanLin}) are of type two and type three. 
These invariants have been studied from various angles: for
instance, combinatorial formul\ae\ for evaluating them have been
derived, and simple bounds in terms of crossing number have been
obtained (see e.g.\ \cite{PolyakViro,Lannes:Degree3,Willerton:Thesis}).  
In this work they are
examined from the novel point of view of the actual values that they take on
knots of small crossing number.  For instance, one can ask how
accurate the known bounds
are, as in Section~1.  When looking at
this question I was led to plot the values of these invariants
which revealed the interesting ``fish'' plots of Section~2: these
pictures form the focus of this note.  Various
ideas which arise from these graphs can be answered for torus knots
and this is done is Section~3.  The final section presents some
problems and further questions.

\section{On $v_2$ and $v_3$.}
The space of additive invariants of type three is two dimensional.
By ``the first two Vassiliev invariants'' is meant the elements of a
basis $\{v_2,v_3\}$ of this space.  The invariants $v_2$ and $v_3$ can
be defined canonically in the following fashion.
The space of additive invariants of type three splits
into the direct sum of type three invariants which do not distinguish
mirror image knots and the type three invariants which differ by a
factor of minus one on mirror image knots.
%   \footnote{I.e.\ this splitting is into the plus one and minus one
%   eigenspaces of the 
% adjoint of the mirror image map.}
Pick the vector in each of these
one dimensional spaces which takes the value one on the positive trefoil.
The one which is invariant under taking mirror images is of type two
and will be denoted $v_2$; the other will be denoted $v_3$.

The invariant $v_2$ has appeared in various guises previously in knot
theory: it is the coefficient of $z^2$ in the Conway polynomial
and its reduction modulo two is the Arf invariant.
Both $v_2$ and $v_3$ can be obtained from the Jones polynomial in the
following fashion.  If $J(q)$
is the Jones polynomial of a knot $K$ and $J^{(n)}(q)$ denotes the
$n^{\text{th}}$ derivative with respective to $q$, then 
  $$v_2(K)=-\tfrac{1}{6}J^{(2)}(1),\qquad
           v_3(K)=-\tfrac{1}{36}\left(J^{(3)}(1)+ 3J^{(2)}(1)\right).$$

Combinatorial formul\ae\ for $v_2$ and $v_3$ can be given in terms of
Gau\ss\ diagram formul\ae --- the reader is referred to
\cite{PolyakViro,Willerton:Thesis}.  
From the combinatorial formul\ae, it is straight forward to obtain simple
bounds for $v_2$ and $v_3$ in terms of the crossing number, $c$, of
the knot, $K$:
namely
  $$|v_2(K)|\le \tfrac{1}{4}c(c-1),\qquad
             |v_3(K)|\le \tfrac{1}{4}c(c-1)(c-2).$$
The first of these bounds was obtained by Lin and Wang \cite{LinWang:IntGeom}
and led Bar-Natan \cite{Bar:Poly} to prove that any type $n$
invariant is bounded by a degree $n$ polynomial in the crossing number
--- this also follows from Stanford's algorithm
\cite{Stanford:ComputingVIs} for calculating
Vassiliev invariants.  The bound for $v_3$ was obtained in
\cite{Willerton:Thesis} by utilizing Domergue and Donato's integration
\cite{Do-Do}
of a type three weight system.

It is natural to ask how sharp these bounds are, and it is this
question that motivated this work.  Stanford has calculated
Vassiliev invariants up to order six for the prime knots up to ten
crossings, the programs and data files of which are available as
\cite{Stanford:Programs}; Thistlethwaite has calculated various
polynomials for knots up to fifteen crossings, these are available in
the {\tt knotscape} program \cite{HosteThistlethwaite:knotscape}.
Using these data, one can compare the bounds on $|v_2|$ and $|v_3|$
given above, with the actual maximum attained for each crossing number
--- this comparison is made in Table \ref{maxmintable}.  It is seen
that in this range of crossing numbers, the bounds are not
particularly tight.

\begin{table}\begin{center}

\begin{tabular}{|l||c|c|c|c|c|c|c|c|c|c|}
\hline
Crossing number&3&4&5&6&7&8&9&10&11&12\\
\hline\hline
Maximum $|v_2|$&1&1&3&2&6&5&10&9&15 &14\\
  Bound  on $|v_2|$& 1.5 &  2 & 5 & 7.5 & 11.5 & 14 &18 & 22.5&27.5 & 33\\
 \hline
 Maximum $|v_3|$&1& 0& 5& 1 & 14 & 10& 30& 25&55 &49 \\
Bound on $|v_3|$& 1.5 & 6   & 15 & 30  & 57.5& 84 & 126 & 180& 247.5 &330\\
\hline
\end{tabular}
\end{center}
\caption{Comparing actual maxima and minima of $|v_2|$ and $|v_3|$ with
the bounds of Section~1.}
\label{maxmintable}
\end{table}

By looking at the raw data, one can see that in this range, for
 odd crossing number
$(2b+1)$, the maxima are achieved precisely by the $(2,2b+1)$-torus
knot, and that this dominates the $v_2$ and $v_3$ of the
$(2b+2)$-crossing knots as well.
Letting $T(p,q)$ be the knot type of
the $(p,q)$-torus knot, Alvarez and Labastida
\cite{AlvarezLabastida:Torus} (see also Section~3 below) give
explicitly for crossing number $c=2b+1$,
$$v_2\left( T\left(2,c\right)\right)=(c^2 -1)/8,\qquad v_3\left(
  T\left(2,c\right)\right)=c(c^2 -1)/24.$$
One could conjecture that
these give bounds on $v_2$ and $v_3$.  After an earlier version of
this paper, Polyak and Viro \cite{PolyakViro:OnTheCassonKnotInvariant}
showed that for a knot with $c$ crossings $v_2\le c^2/8$.

\section{Plots for knots with up to fourteen crossings.}
\begin{figure}
\begin{minipage}[b]{.3\linewidth}
  %\hskip -.1in
     \epsfig{file=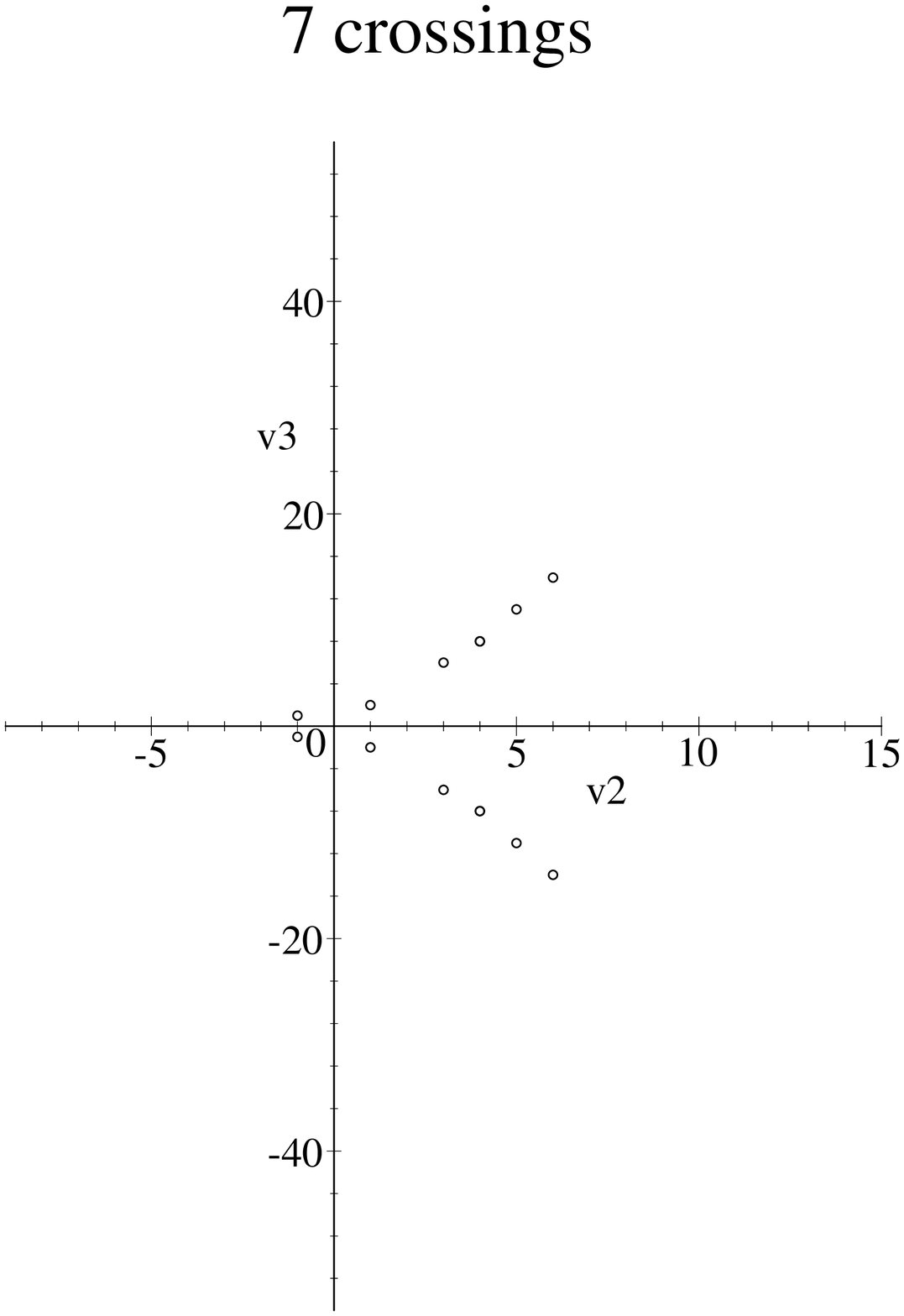,width=1.2\linewidth}
\end{minipage}\hfill
\begin{minipage}[b]{.3\linewidth}
  %\hskip -.1in
     \epsfig{file=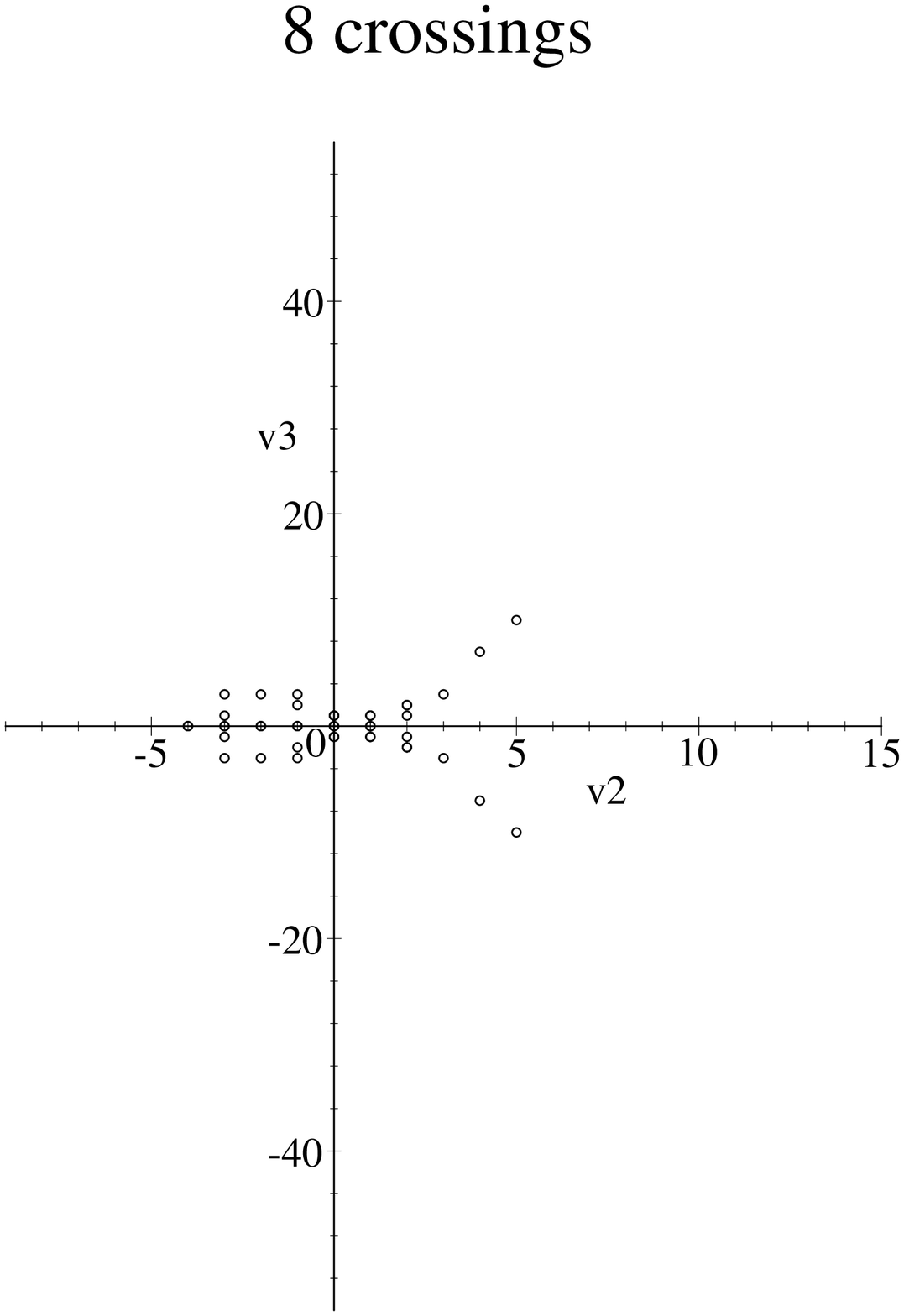,width=1.2\linewidth}
\end{minipage}\hfill
\begin{minipage}[b]{.3\linewidth}
  \hskip -.1in
     \epsfig{file=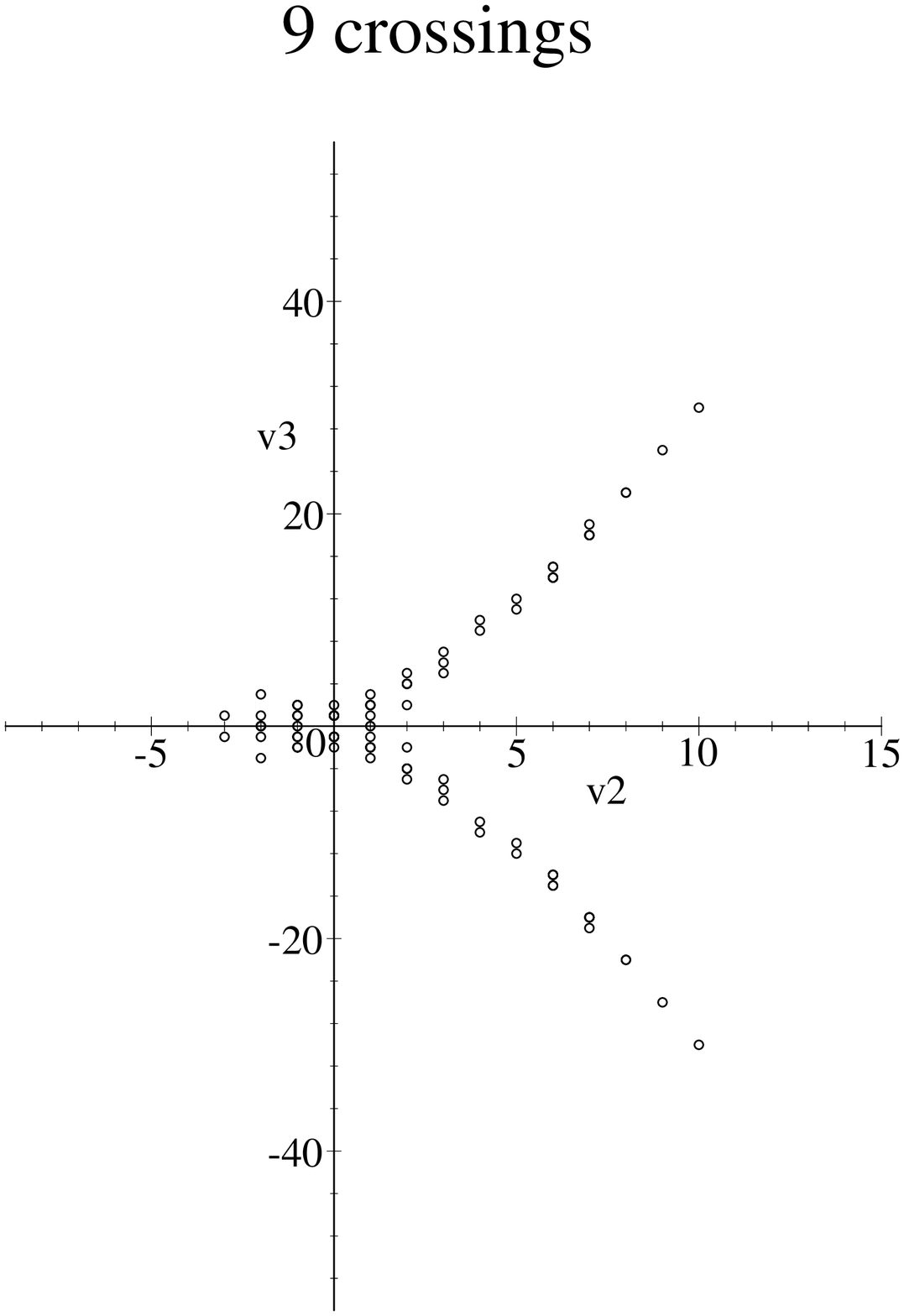,width=1.2\linewidth}
\end{minipage}\hfill\vskip -.1in
\begin{minipage}[b]{.3\linewidth}
  \hskip -.1in
     \epsfig{file=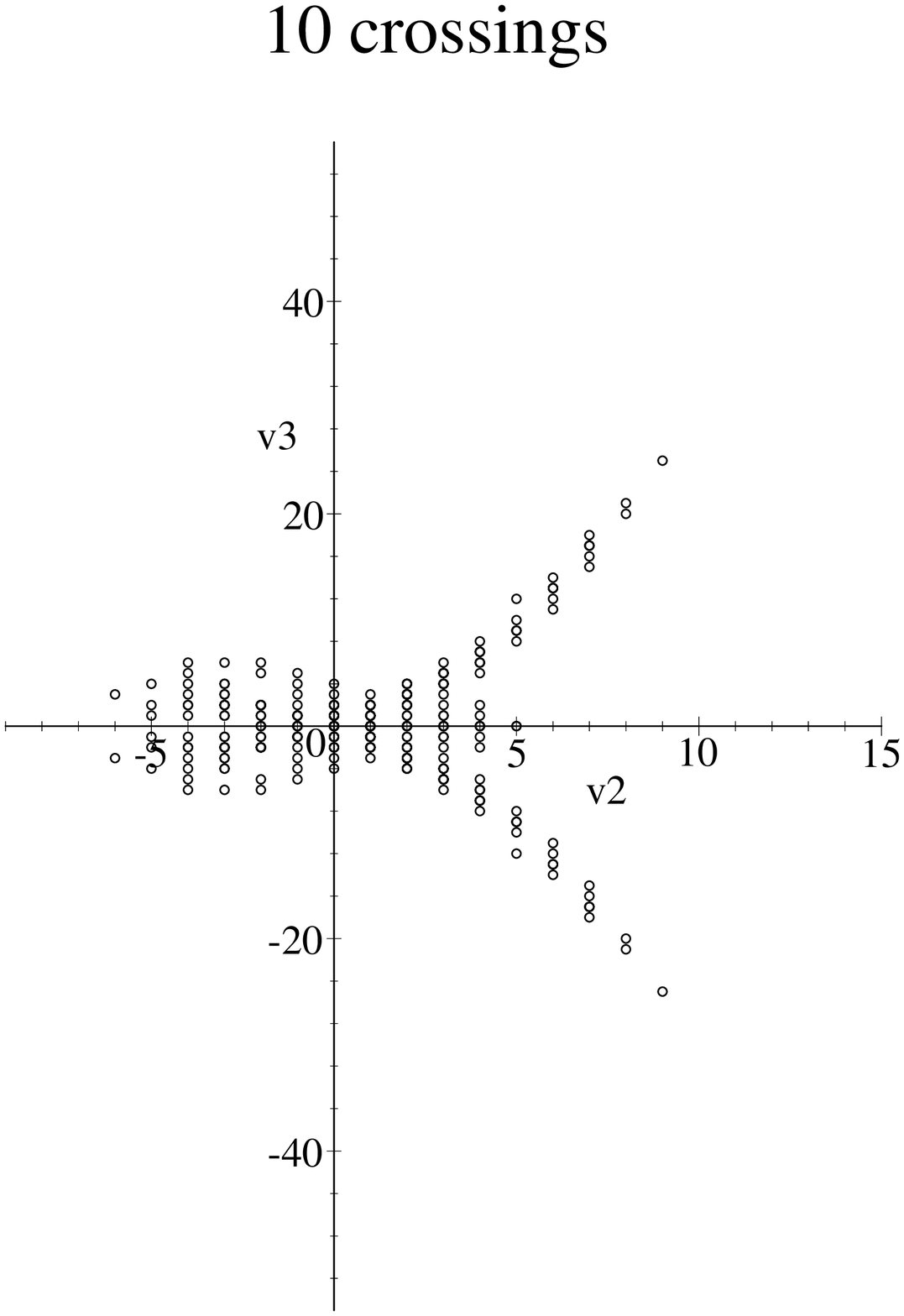,width=1.2\linewidth}
       
\end{minipage}\hfill
\begin{minipage}[b]{.3\linewidth}
  \hskip -.1in
     \epsfig{file=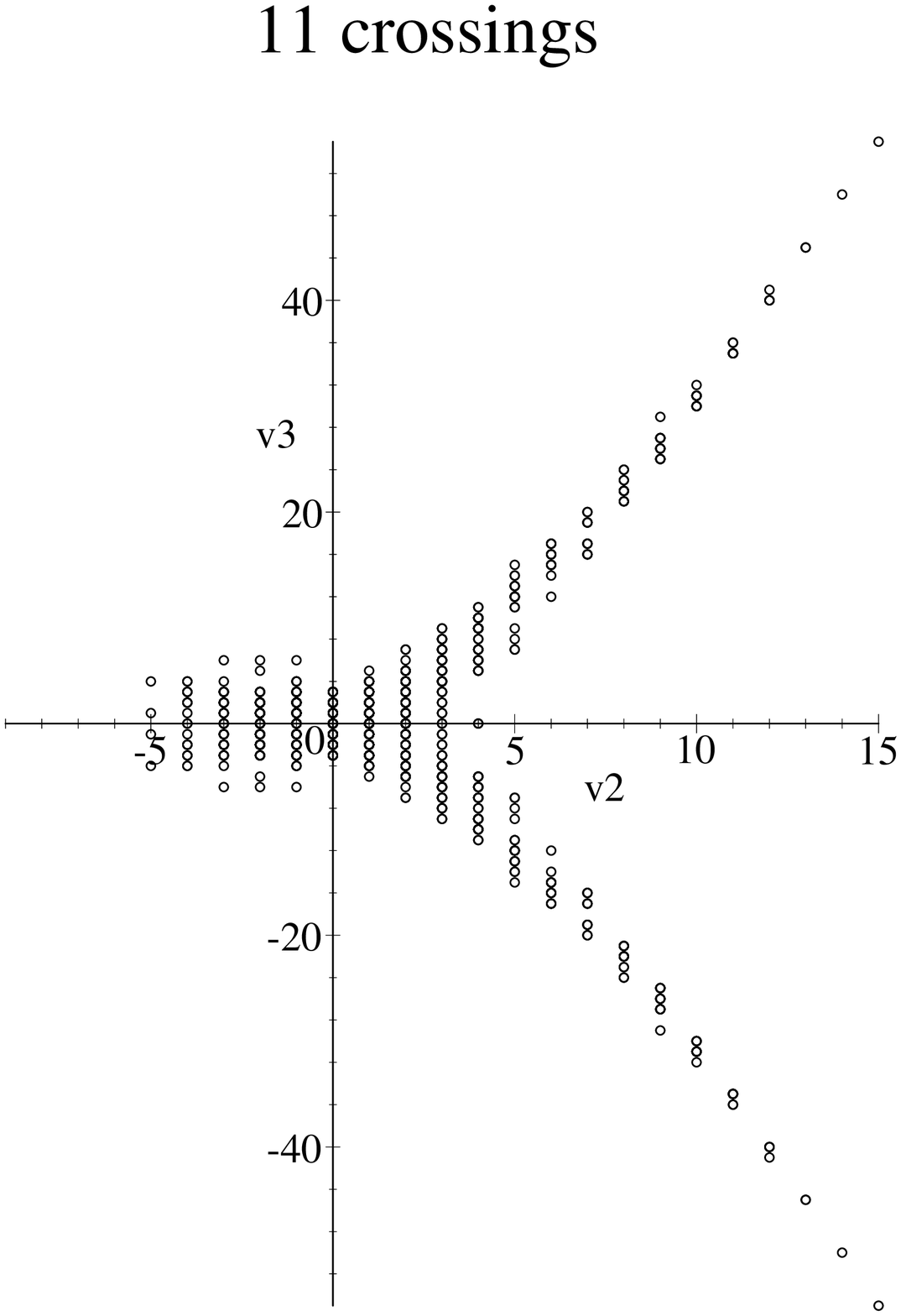,width=1.2\linewidth}
\end{minipage}\hfill
\begin{minipage}[b]{.3\linewidth}
  \hskip -.1in
     \epsfig{file=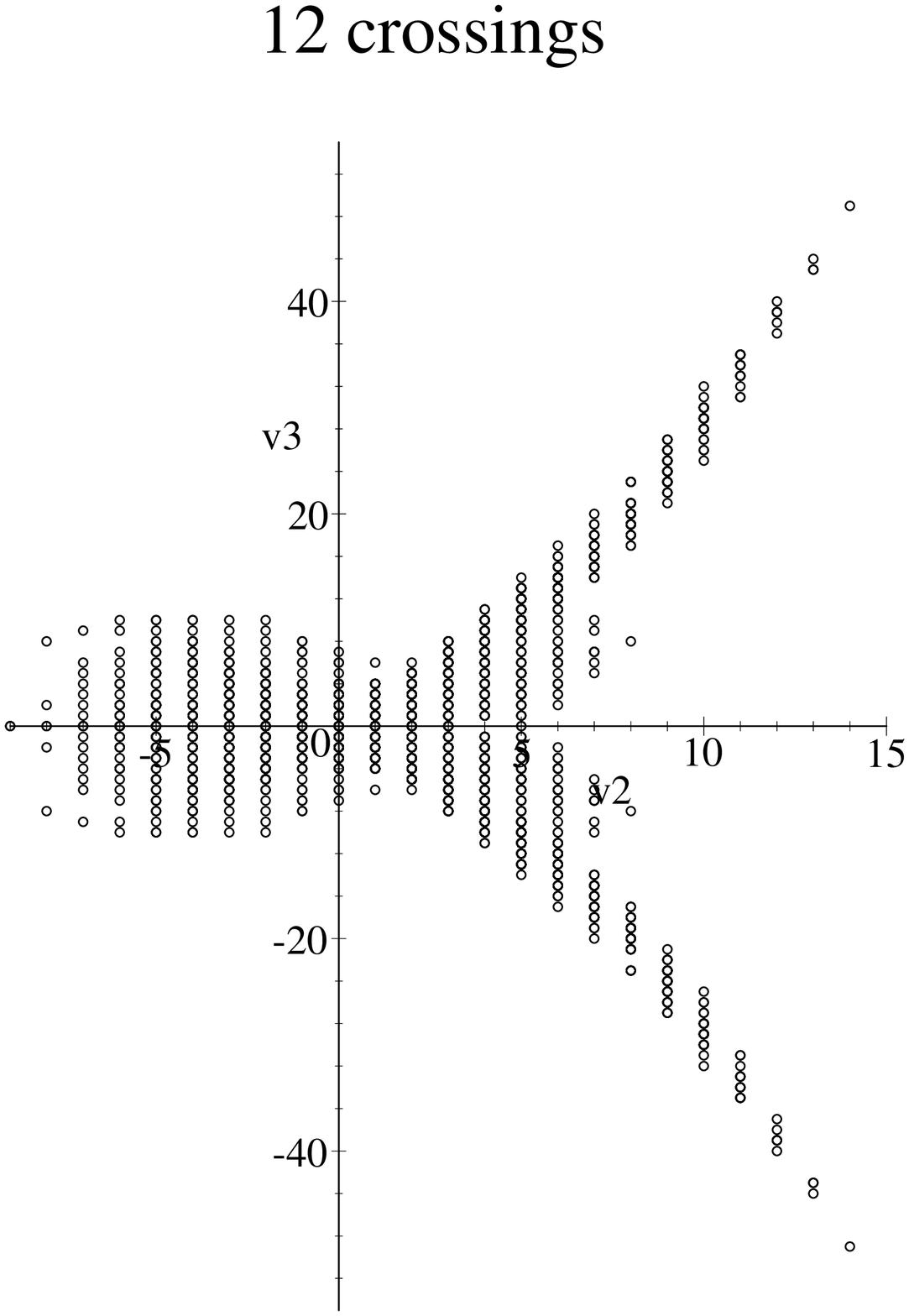,width=1.2\linewidth}
\end{minipage}\hfill
\caption{Plots by crossing number of $v_2$ and $v_3$ for the prime
knots up to twelve crossings.}
\label{v2v3plotsto12}
\end{figure}
\begin{figure}
\begin{minipage}[b]{.45\linewidth}
\epsfig{file=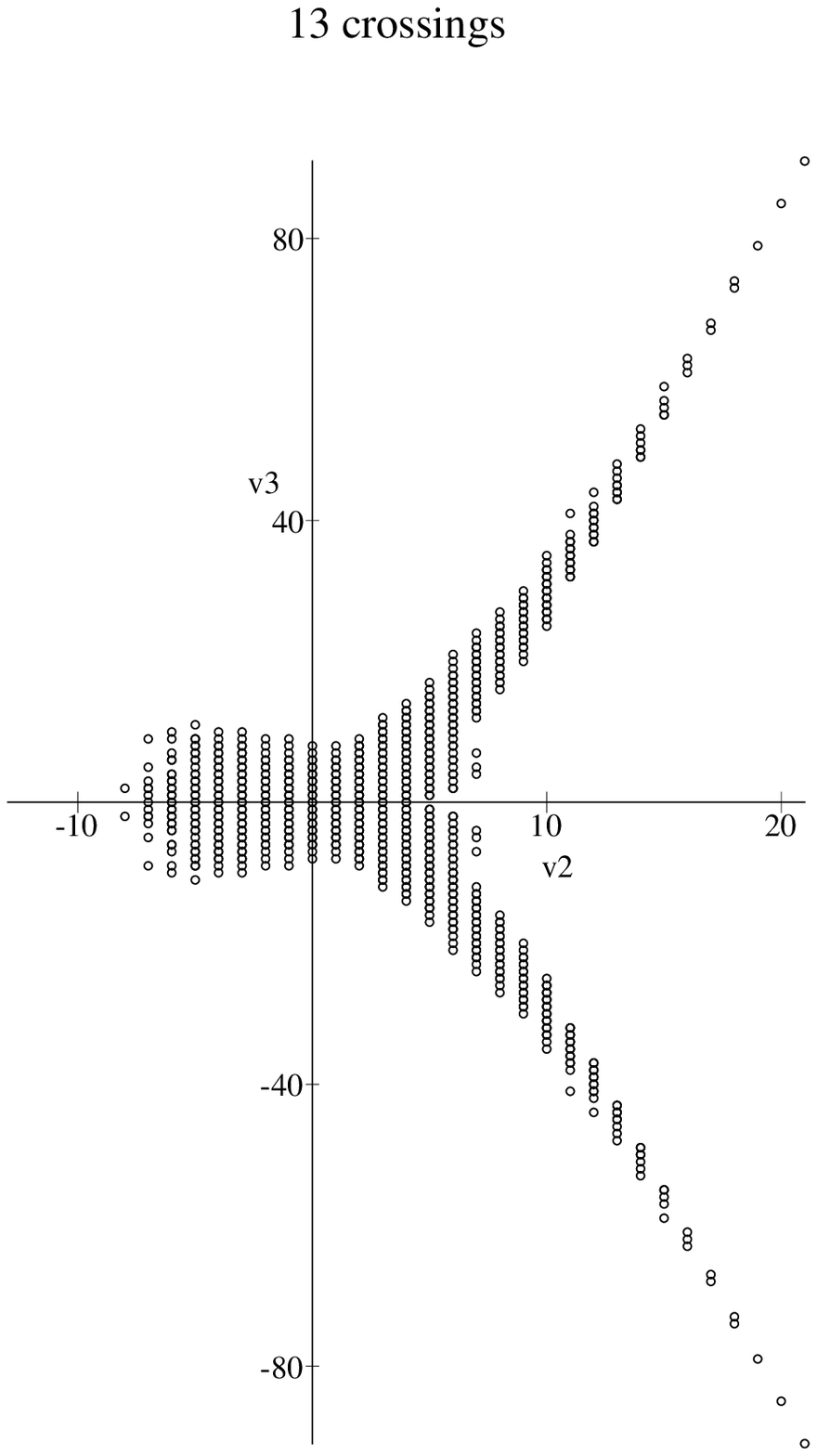,width=1.2\linewidth}
\end{minipage}
\begin{minipage}[b]{.45\linewidth}
\epsfig{file=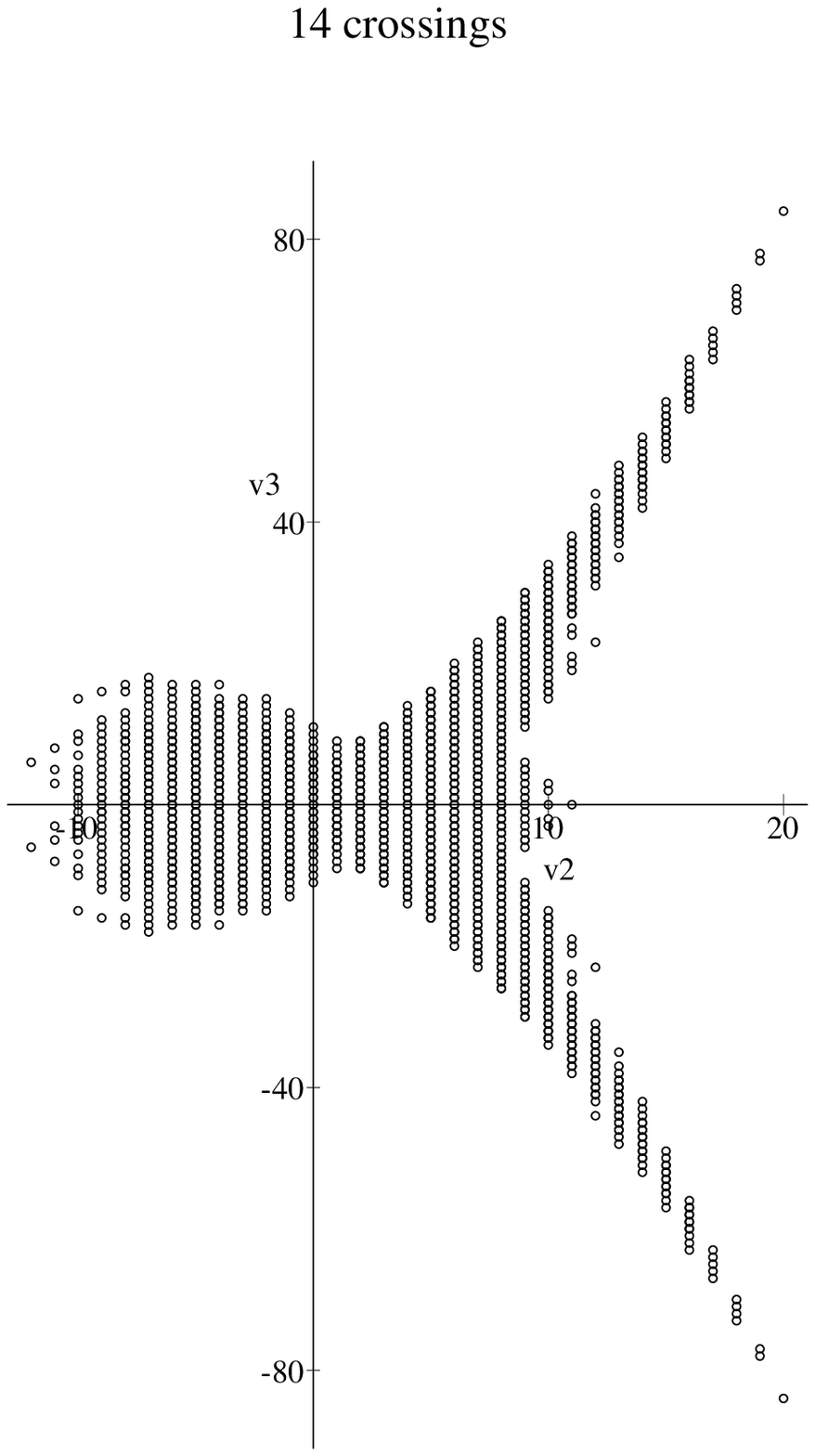,width=1.2\linewidth}
\end{minipage}
\vskip-.5in
\caption{Plots by crossing number of $v_2$ and $v_3$ for the prime
knots with thirteen and fourteen crossings.  
}
\label{v2v3plots13and14}
\end{figure}

Having stared at the raw data of Stanford for sufficiently long to
start noticing patterns, I was led to plot $v_2$ against $v_3$ for
knots of each crossing number up to  crossing number fourteen.  These
plots are contained in Figure~\ref{v2v3plotsto12} and
Figure~\ref{v2v3plots13and14}.  The symmetry in the $v_2$-axis is
expected, as this is just the effect of taking the mirror image of the
knots.  The ``fish'' shape of these plots is not expected! 
This shape suggests
some bound of the form
  $$\text{cubic in\ } v_2(K) \le (v_3(K))^2 \le \text{another cubic 
      in\ } v_2(K).$$
Such bounds, independent of crossing number do in fact exist for torus
knots, as will be seen below.  However, this 
cannot be the case in
general (unless the bounds depend on the crossing number): two reasons
are as follows.   

\begin{figure}[t]
$$\vpic{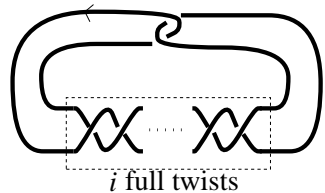}$$
\caption{The $i$th twisted Whitehead double of the unknot, $\Wh (i)$.
For $i$ negative, $i$ full twists means $-i$ negative twists.}
\label{Whiteheadfig}
\end{figure}

Firstly, consider the sequence of Whitehead doubles of the unknot,
$\{\Wh (i)\}_{i\in \Z}$ (see Figure~\ref{Whiteheadfig}).  
Table~\ref{Whiteheadvalues} gives the value
of $v_2$ and $v_3$ on these for a range of $i$.  It follows from the
theorem of Dean \cite{Dean} and Trapp \cite{Trapp} on twist sequences
that a type $n$ invariant evaluated on the Whitehead doubles is a
polynomial in $i$ of degree at most%
\footnote{By an observation of Lin in this case it must be of degree
at most $n-1$.}
$n$.  A glance at Table~\ref{Whiteheadvalues}
suffices to deduce that $v_2\left( \Wh \left(i\right)\right)=i$ and
$v_3\left( \Wh\left(i\right)\right)=\tfrac{1}{2}i(i+1)$.  Thus there
is a sequence of knots (all except the unknot having unknotting number
equal to one) which maps into the $(v_2,v_3)$-plane as a nice
quadratic.  This contradicts any bounds of the above form.
\begin{table}[t]          
$$\begin{array}{|c||c|c|c|c|c|c|c|c|}
\hline
i & -3 &-2 & -1& 0 & 1 & 2 & 3& 4\\
\hline
\Wh(i) &8_1&6_1&4_1&0_1&3_1&5_2&7_2&9_2\\
\hline
v_2\left(\Wh(i)\right) & -3 &-2 & -1& 0 & 1 & 2 & 3& 4\\
v_3\left(\Wh(i)\right) & 3 &1 &0 &0&1&3&6&10\\
\hline
\end{array}$$
\caption{The values of $v_2$ and $v_3$ on the twisted Whitehead doubles
of the unknot.  The knot notation, e.g.\ $3_1$, refers to
Alexander-Briggs notation (see \cite{BurdeZieschang}).}
\label{Whiteheadvalues}
\end{table}

Secondly, for any $(a,b)\in \Z^2$ one can obtain a prime (alternating) knot
with $(v_2,v_3)$ equal to $(a,b)$ in the following manner: connect
sum suitably many positive and negative trefoil knots (with
$(v_2,v_3)=(1,\pm 1)$) and figure eight knots (with $(v_2,v_3)=(-1,0)$),
to obtain a composite knot 
%\nopagebreak[0] 
with $(v_2,v_3)=(a,b)$, then Stanford
\cite{Stanford:BraidCommutators} gives a method for constructing a prime
knot with the same $v_2$ and $v_3$.  

There does appear to be a qualitative difference between the pictures
for odd and even crossing numbers in Figures~\ref{v2v3plotsto12}
and~\ref{v2v3plots13and14}.  The even
crossing number ones seem to be more concentrated in the `body' of the `fish'
and the odd ones more in the `tail'.  Note that for each odd crossing
number, $c$, there is the $(2,c)$-torus knot and the Whitehead double
$\Wh (\left(c-1\right)/2)$ with a $(v_2,v_3)$ of 
$\left(\left(c-1\right)/2,\left(c^2-1\right)/8\right)$; and for
even crossing number, $c$, there is the Whitehead double $\Wh (1-c/2)$
with a $(v_2,v_3)$ of $\left(1-c/2,\left(c-2\right)c/8\right)$.  
Also for up to twelve
crossings the amphicheiral knots
--- that is those equivalent to their mirror image, and hence with
$v_3=0$ --- all have even crossing number, but this is not true in
general, as the fifteen crossing knot $15_{224980}$ is amphicheiral.

\section{Torus knots.}
\begin{figure}
\hfill(i)\begin{minipage}{.45\linewidth}\epsfig{file=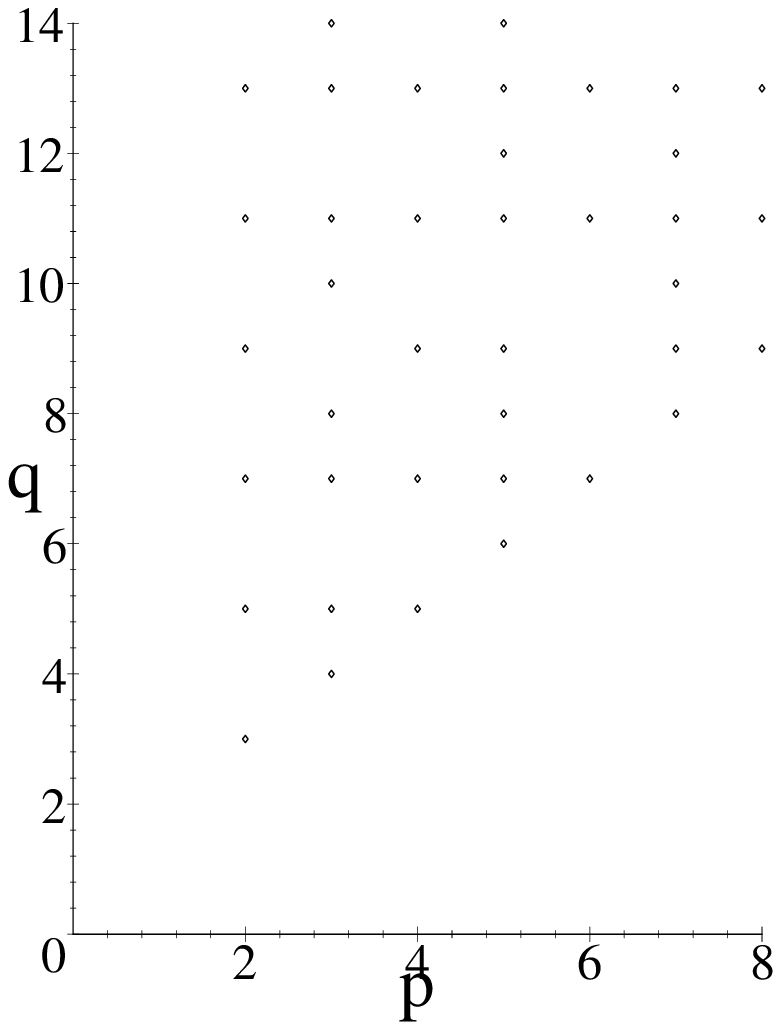,width=\linewidth}
\end{minipage} \kern-20pt
$\longrightarrow$ \kern-10pt
\begin{minipage}{.45\linewidth}
 \epsfig{file=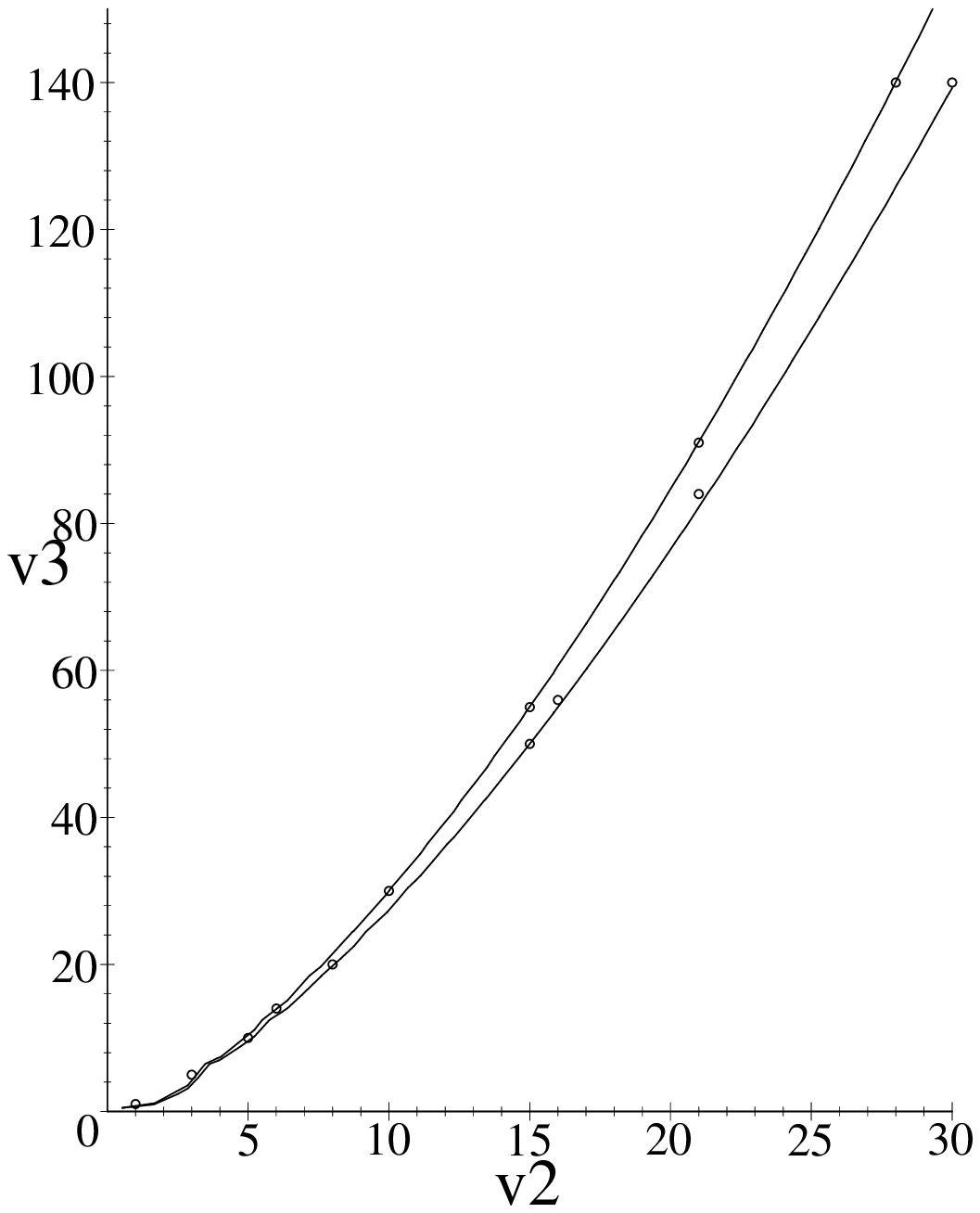,width=\linewidth}
\end{minipage} \hfill
\\ \vskip -.2in
(ii)\begin{minipage}{.45\linewidth}
\epsfig{file=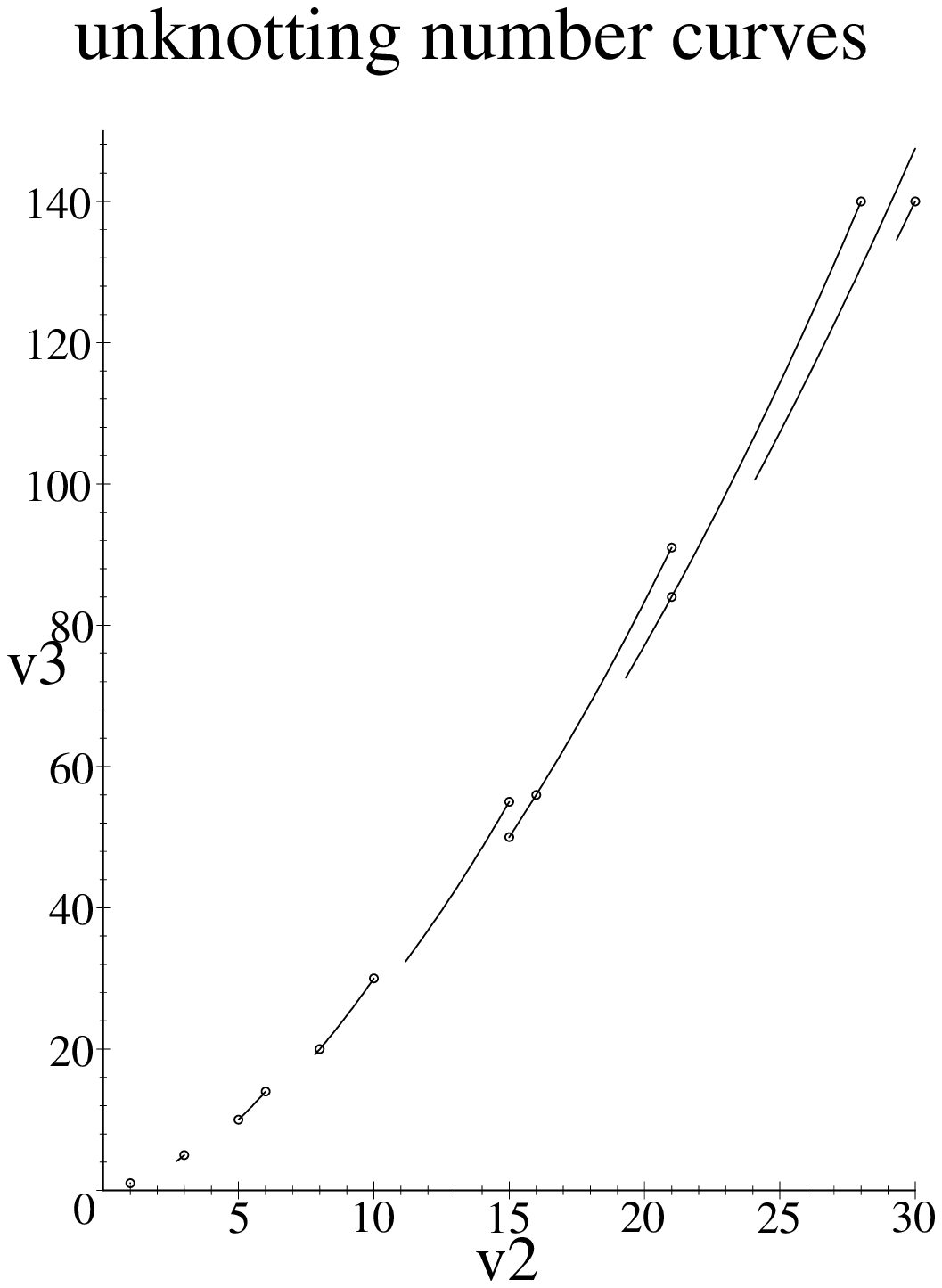,width=\linewidth}\end{minipage}
\hfill
(iii)\begin{minipage}{.45\linewidth}
\epsfig{file=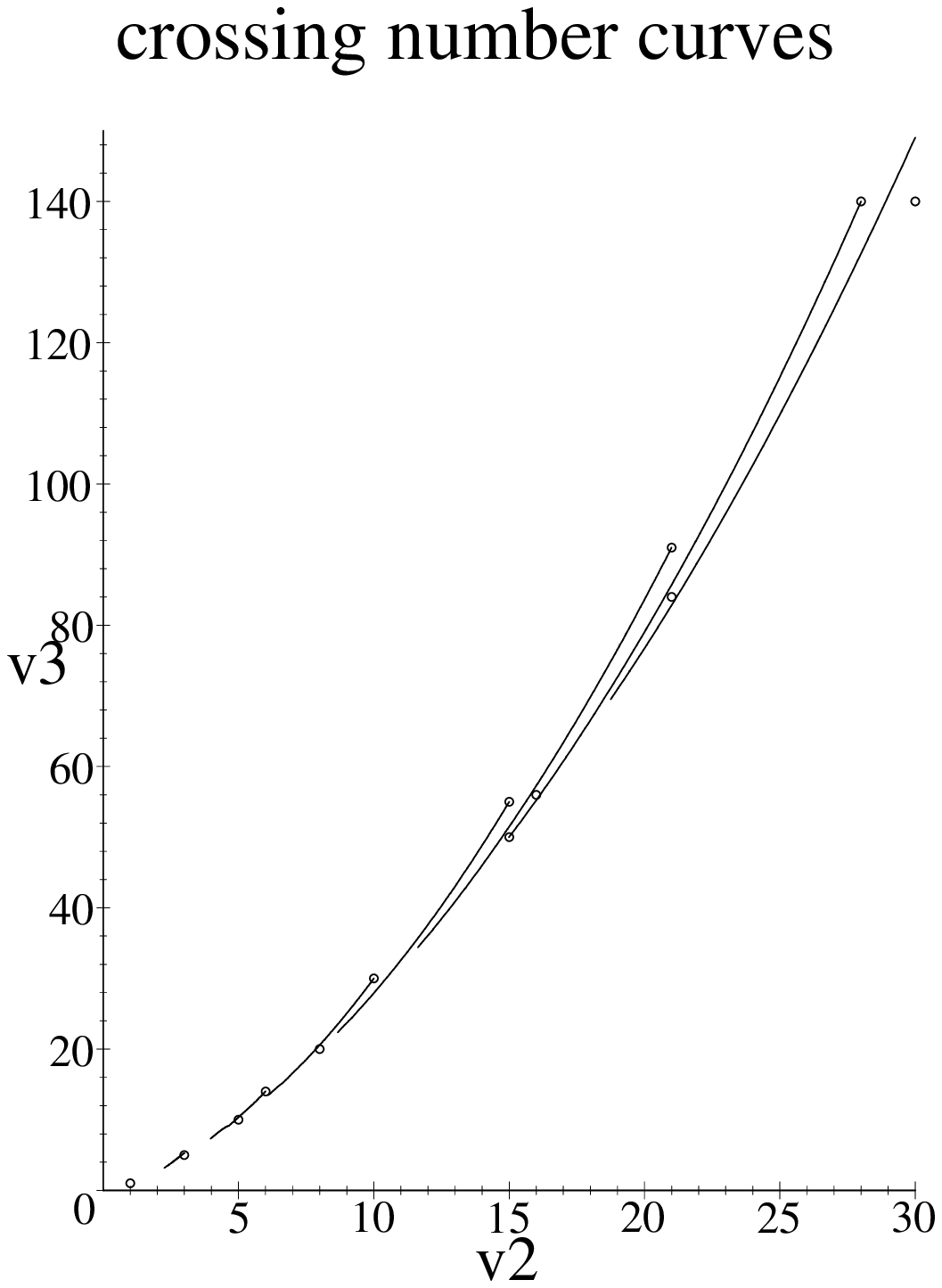,width=\linewidth}\end{minipage}
%$$\text{(ii)}\vpic{v_torusunknotting}
%      \text{(iii)}\vpic{v_toruscrossing}$$
\caption{Torus knots in the $(v_2,v_3)$-plane: (i) mapping torus knots
from the $(p,q)$-plane into the region of the $(v_2,v_3)$-plane given by
Propositions \ref{cubicbounds1} and \ref{cubicbounds2};  (ii) torus
unknotting number curves for $u=1,\dots,9$ (see Section~3.2); (iii)
torus crossing number curves for $c=3,5,\dots,17$ (see Section~3.3).}
\label{torusgraphs}
\end{figure}

The purpose of this section is to show that the torus knots map into
the $(v_2,v_3)$-plane in a nice manner.  In particular they satisfy
cubic bounds of the form described above, implying that they lie on
the tails of the fish; further, torus knots of the same unknotting
number, or crossing number, lie on nice curves in the
$(v_2,v_3)$-plane.  The results of this section are summarized
diagrammatically in Figure~\ref{torusgraphs}.
\renewcommand{\topfraction}{.99}

For $p$ and $q$ coprime let $T(p,q)$ be the knot type of the
$(p,q)$-torus knot.  Then $T(p,q)$ is the
unknot if and only if $p$ or $q$ is $\pm 1$, and for $T(p,q)$
nontrivial, $T(p,q)$ is the same
knot as $T(p',q')$ if and only if $(p',q')$ equals one the following
$(p,q)$, $(q,p)$, $(-p,-q)$, or $(-q,-p)$.  Further $T(p,-q)$ is the
mirror image of $T(p,q)$.  See \cite{BurdeZieschang}.

  The key to this section is the following pair of
formul\ae\ 
of Alvarez and Labastida \cite{AlvarezLabastida:Torus}.
  $$v_2(T(p,q))=\tfrac{1}{24}(p^2-1)(q^2-1),\qquad
                  v_3(T(p,q))=\tfrac{1}{144}pq(p^2-1)(q^2-1).$$
Note that these have the required properties under the symmetries of $p$
and $q$ mentioned above, and that these are integer valued on torus
knots (i.e.\ when $p$ and $q$ are coprime).  Also $T\mapsto
(v_2(T),v_3(T))$ is injective for torus knots, that is to say torus
knots are determined by their $(v_2,v_3)$.

\subsection{Cubic bounds}
With the above formul\ae\ of Alvarez and Labastida it is
straightforward to prove bounds, for 
torus knots, of the
form suggested in the last section.
\begin{prop}
If $T$ is a torus knot then 
  $$\frac{2}{3} v_2(T)^3 +\frac{1}{3}v_2(T)^2 \le v_3(T)^2
                      \le \frac{8}{9} v_2(T)^3 +\frac{1}{9}v_2(T)^2.$$
%\nopagebreak[0]%
Further, the right hand bound is tight in the sense that there exist
torus knots with \nopagebreak[0] arbitrarily large $v_2$ and $v_3$ such that equality
holds.
\label{cubicbounds1}
\end{prop}
\begin{proof}
Suppose that $T$ is a $(p,q)$-torus knot then
\begin{align*}
 v_3(T)^2-\frac{2}{3}v_2(T)^3&= \left(\frac{1}{144} pq(p^2-1)(q^2-1)\right)^2
                 -\frac{2}{3}\left(\frac{1}{24}(p^2-1)(q^2-1)\right)^3\\
          &= \frac{1}{12^4}(p^2-1)^2(q^2-1)^2 [p^2+q^2 -1]         \\
          &\ge \frac{1}{12^3}(p^2-1)^2(q^2-1)^2 
                  \qquad\qquad\text{\ as\ }   p^2+q^2 \ge 13          \\
          &= \frac{1}{3}v_2(T)^2,
\end{align*}
hence the first inequality (with equality only in the case of torus knots).

For the second,
\begin{align*}
\frac{8}{9}v_2(T)^3-v_3(T)^2&=\frac{8}{9}\left(\frac{1}{24}(p^2-1)(q^2-1)\right)^3
           - \left(\frac{1}{144} pq(p^2-1)(q^2-1)\right)^2   \\
        &=  \frac{1}{4.27}\left[\frac{1}{24}(p^2-1)(q^2-1)\right ] ^2
           \left\{ 4(p^2-1)(q^2-1) - 3 p^2 q^2\right\}                  \\
        &= \frac{1}{4.27}v_2(T)^2 \left\{(p^2-4)(q^2-4)-12\right\}      \\
        &\ge \frac{1}{4.27}v_2(T)^2 \{-12\} = -\frac{1}{9} v_2(T)^2,
\end{align*}
and note that equality occurs precisely when $T$ is a $(2,q)$-torus knot.
\end{proof}

Although the left hand bound has the correct asymptotic behaviour, for a
tight bound a different form of cubic is required.

\begin{prop}
For a torus knot $T$,
  $$\frac{2}{3}v_2(T)^3 +\frac{1}{3}v_2(T)v_3(T) \le v_3(T)^2,$$
and this bound is tight in the sense of the previous proposition.
\label{cubicbounds2}
\end{prop}
\begin{proof}
Using the notation of the previous proof,
\begin{align*}
v_3(T)^2 -\frac{2}{3}v_2(T)^3 -\frac{1}{3}v_2(T) v_3(T) &= \frac{1}{36.24^2}
                    (p^2-1)^2 (q^2-1)^2 \left ( (p-q)^2 -1\right) \\
%%         &= \frac{1}{36} v_2^2 \left ( (p-q)^2 -1\right)\\
         &\ge 0,
\end{align*}
with equality if and only if $T$ is a $(p,p+1)$ torus knot.
\end{proof}

Given that half the torus knots (those with positive $v_3$) can be
thought of as lying in the region $q>p>0$ in the $(p,q)$-plane, these
bounds are not surprising.  Graphically this can be seen in Figure
\ref{torusgraphs}.

\subsection{Torus knots and unknotting number}
By Kronheimer and Mrowka's \cite{KronheimerMrowka} positive solution to the 
Milnor conjecture the following formula is known for the unknotting
number, $u$, of torus knots:
  $$u(T(p,q))=\frac{1}{2}\left( |p|-1 \right) \left( |q|-1 \right).$$
As a consequence, the following easily verifiable relationship is
obtained:
\begin{prop}
For a torus knot T,
  $$v_2(T)^2+\frac{1}{6}u(T)(u(T)-1)v_2(T)=u(T)|v_3(T)|,$$
and given $v_2(T)$ and $v_3(T)$ then $u(T)$ is the smaller of the two
roots.      \qed
\end{prop}
\noindent So for a fixed unknotting number, 
the torus knots lie on a quadratic in
the $(v_2,v_3)$-plane (c.f.\ the Whitehead knots in Section~2).  This is
pictured in Figure \ref{torusgraphs}.  The segments of curves shown
were chosen by the following proposition.

\begin{prop}
For a torus knot $T$,
  $$\frac{1}{2}u(T)(u(T)+1) \ge v_2(T) \ge 
         \frac{1}{6}u(T)\left(u(T)+\sqrt{8u(T)+1} + 2\right),$$
and both bounds are tight.
\end{prop}
\begin{proof}
If $T$ is a $(p,q)$-torus knot, then a minimal amount of manipulation
gives
\begin{align*}
\frac{1}{2}u(T)(u(T)+1) -v_2(T) &= \frac{1}{12}(|p|-1)(|q|-1)  
                                                   (|p|-2)(|q|-2)\\
                &\ge 0,
\end{align*}
with equality if and only if $T$ is a $(2,q)$-torus knot.

For the right hand bound, firstly, let $a$ and $b$ be distinct positive
integers, then $(a-b)^2\ge
1$, so $(a+b)^2\ge 4ab+1$ and thus $a+b\ge \sqrt{4ab+1}$, with equality
precisely when $a$ and $b$ differ by one.

Now for $T$ a $(p,q)$-torus knot,
\begin{align*}
v_2(T)-\frac{1}{6}u(T)&\left(u(T)+\sqrt{8u(T)+1} + 2\right)\\
   &= \frac{1}{12}(|p|-1)(|q|-1)
       \left\{|p|+|q|-2 -\sqrt{4(|p|-1)(|q|-1) +1} \right\} \\
   &\ge 0, 
\end{align*}
by putting $a=|p|-1,\ b=|q|-1$ in the above paragraph.  Note that
equality occurs precisely when $T$ is a $(p,p+1)$-torus knot.
\end{proof}
Weakening the right hand bound to $v_2\ge
\frac{1}{6}u(T)\left(u(T)+5\right)$ and inverting the inequalities reveals the
following corollary.
\begin{cor}
For a torus knot $T$,
  $$\sqrt{1+8v_2(T)} -1 \le 2u(T) \le \sqrt{24v_2(T) +25} -5,$$
and the left hand bound is tight (in the sense of Proposition~\ref{cubicbounds1}).   \qed
\end{cor}

\subsection{Torus knots and crossing number.}
By the work of Murasugi \cite{Murasugi}, a similar formula is known for
the crossing number, $c$, of torus knots:
  $$c(T(p,q))=|q|(|p|-1), \qquad \text{when\ }|p|<|q|.$$
This leads to the following relation;
\begin{prop}
If $T$ is a torus knot, and 
     $\rho(T)=\left| \frac{6v_3(T)}{v_2(T)}\right|$
then
  $$24v_2(T)(c(T)-\rho(T))^2=c(T)\left(\left(c(T)-\rho(T)\right)^2 -1 \right)
                            \left(2\rho(T) -c(T)\right),$$
and
  $$c(T)=\rho(T)-\frac{1}{2}\left( \sqrt{(\rho(T)-1)^2 -24v_2(T)}
                                  + \sqrt{(\rho(T)+1)^2 -24v_2(T)}\right).$$
\end{prop}
\begin{proof}
This is easily verified; note that if $T$ is a $(p,q)$-torus knot then
$\rho(T)=|pq|$ and $c(T)-\rho(T)=|q|$.
\end{proof}
This isn't as nice a relationship as with the unknotting number: for a fixed
crossing number the relationship is a not particularly nice quartic
between $v_2$ and $v_3$.  However, the crossing number curves can
still be graphed, as in Figure \ref{torusgraphs} --- the length of arc
segments plotted there being determined by the following proposition.
\begin{prop}
For a torus knot $T$,
  $$\frac{1}{8}\left(c(T)^2-1\right)\ge v_2(T)\ge \frac{1}{24}
  c(T)\left( c(T)+1+2\sqrt{c(T)+1}\right),$$
and these bounds are tight (in the sense of Proposition~\ref{cubicbounds1}).
\end{prop}
\begin{proof}
Suppose that $T$ is a $(p,q)$-torus knot with $q>p>0$ --- this just avoids
excessive modulus signs in the calculation ---  then for the left hand
bound,
\begin{align*}
\frac{1}{8}\left(c(T)^2-1\right) -v_2(T)
 &= \frac{1}{24}\left\{ 3\left(\left[ q(p-1)\right]^2 -1 \right)
     -(p^2-1)(q^2-1)\right\}          \\
 &=  \frac{1}{24}\left\{ 2q^2 p^2-6q^2p + 4q^2 +p^2 -4\right\}  \\
 &= \frac{1}{24}(p-2)\left\{(2q^2+1)(p-1)+3 \right\}\\
 &\ge 0,
\end{align*}
and equality occurs precisely when $T$ is a $(2,q)$-torus knot.

For the right hand bound,
\begin{align*}
24 v_2(T) -c(T) &\left( c(T)+1  +2 \sqrt{c(T)+1} \right )    \\
 &= (p^2-1)(q^2-1)-q(p-1)\left(q(p-1)+1+2\sqrt{q(p-1)+1} \right) \\
 &= (p-1)\left\{ 2q^2-q-1-p-2q\sqrt{qp-q+1}\right\},
\end{align*}
and claim that this is non-negative and is zero precisely when $q=p+1$.

To prove the claim, note
  $$(q-1)^2=q(q-1)-q-1 \ge qp-q-1$$
as $q-p-1\ge 0$, and so also
  $$(q-1)^2+\frac{2(q-1)\left( q-p-1\right)}{2q}
    +\left[\!\frac{q-p-1}{2q}\!\right]^2 \geq qp-q-1>0,$$
thus, by taking square roots,   
  $$(q-1)+\frac{q-p-1}{2q} \ge \sqrt{qp-q-1},$$
from which the claim follows on multiplying through by $2q$.
\end{proof}

Weakening the right hand bound to $v_2\ge \frac{1}{24}c(c+5)$ and
inverting, gives
\begin{cor}
For a torus knot T
 $$\frac{1}{24}\left( \sqrt{25+96v_2(T)}-5\right)\ge
          c(T)\ge 2\sqrt{8v_2(T)+1},$$
and the right hand bound is tight in the previous sense.
\qed
\end{cor}

\section{Problems and further questions.}

\begin{problem}
The invariants $v_2$ and $v_3$ form, in some sense, a canonical basis
for the space of invariants of degree three.  
Can one find similarly, for example, 
a canonical basis for the space of type four, even,
additive invariants?  Can canonical bases be found for higher order
invariants? 
\end{problem}
\begin{problem}
Does the fish pattern persist in the graphs of knots with higher crossing number?
\end{problem}

\begin{problem}
Is there some qualitative distinction between knots with odd and even
crossing number which explains the perceived difference in the fish?
\end{problem}

\begin{problem}
Is there any relationship with unknotting number?  Note that the
$n$-fold connect sum of $8_{14}$ has unknotting number $n$ and
$(v_2,v_3)=(0,0)$ (this was pointed out to me by Stoimenov).  

\end{problem}

\begin{problem}
For a knot $K$ with $(6|v_3(K)|-|v_2(K)|)^2\ge 24 v_2(K)^3$, let
$\rho(K)=6|v_3(k)/v_2(K)|$ and then define the {\em pseudo-unknotting
number}, $\tilde u(K)$, and the {\em pseudo-uncrossing
number}, $\tilde c(K)$, by
\begin{align*}
\tilde u(K)&:=\half\left(1+\rho(K)-\sqrt{(1+\rho(K))^2 -24 v_2(K)}\right);\\
\tilde c(K)&:=\rho-\half\left(\sqrt{(1+\rho(K))^2 -24 v_2(K)}
                        +\sqrt{(1-\rho(K))^2 -24 v_2(K)} \right).
\end{align*}
For torus knots, the pseudo-unknotting and pseudo-crossing numbers
coincide with the usual unknotting and crossing numbers.  Do they have
any meaning for other knots?  Does the necessary bound for $K$
have any topological interpretation?

As an example, consider the Whitehead knots $\Wh (i)$, for $i>0$ these
all have unknotting number equal to one; in this case $\tilde
u(\Wh(1))=1$, and $\tilde u(\Wh(i))\rightarrow 2$ as $i\rightarrow \infty$. 
\end{problem}
\section*{Acknowledgements}
This work formed part of the author's University of Edinburgh PhD
Thesis, and was supported by an EPSRC studentship.  Thanks to
Oliver Dasbach and Alexander Stoimenow for useful comments.  The
diagrams and numerical experiments were done using Maple.

\bibliographystyle{amsplain}
\bibliography{Vas}
\end{document}